\documentclass[12pt, a4paper]{article}
\usepackage[utf8]{inputenc}
\usepackage[a4paper,left=1.5cm,right=1.5cm,top=2cm,bottom=2cm]{geometry}
\usepackage{authblk}
\usepackage{graphicx}
\graphicspath{ {./images/} }
\usepackage{indentfirst}
\usepackage{natbib}
\usepackage[colorlinks=true,citecolor=blue,urlcolor=blue]{hyperref}
\usepackage{epstopdf}
\usepackage{amssymb}
\usepackage{amsmath}
\usepackage{longtable}
\usepackage{multirow}

\newtheorem{definition}{Definition}

\newtheorem{example}{Example}

\newtheorem{procedure}{Procedure}
\newtheorem{hypothesis}{Hypothesis}

\title{\vspace {-1cm} \large \textbf{Mitigating the rank reversal probability in the entropy-weight TODIM method through an expanded evaluation matrix: A case study on energy storage technology selection}}

\date{}
\author[1]{Lisheng Jiang}
\author[1]{Tianyu Zhang}
\author[1,3]{Shiyu Yan\thanks{Corresponding author. Email: syyan@cdut.edu.cn.}}
\author[2]{Ran Fang}

\affil[1]{{\normalsize College of Management Science, Chengdu University of Technology, Chengdu 610059, China}}
\affil[2]{\normalsize School of Architecture and Civil Engineering, Chengdu University, Chengdu, 610106, China}
\affil[3]{{\normalsize Energy and Environment Carbon Neutrality Innovation Research Center, Chengdu University of Technology, Chengdu 610059, China}}

\begin{document}
\maketitle

\begin{abstract}
The TODIM method (an acronym in Portuguese for interactive and multiple criteria decision-making) with entropy weights is influenced by rank reversal, a phenomenon where the order of two alternatives changes following the addition of another alternative. Research on rank reversal has predominantly focused on single decision-making methods. To the best of our knowledge, the reduction of rank reversal probability in hybrid methods, such as the entropy-weight TODIM method, remains an unresolved challenge. To address this, this paper introduces the expanded evaluation matrix, which incorporates virtual alternatives, to reduce the probability of rank reversal in the entropy-weight TODIM method. A simulation study is conducted to assess the effectiveness of the expanded evaluation matrix in mitigating rank reversal. The results demonstrate that the expanded evaluation matrix significantly reduces the rank reversal probability. A case study on selecting energy storage technology showcases the potential real-world applications of the expanded evaluation matrix. The reliability of the expanded evaluation matrix is further validated through sensitivity and comparative analyses. Given the simplicity and ease of implementation of the expanded evaluation matrix, it can be readily adapted to other decision-making methods and holds substantial potential for broad application.
\end{abstract}


\newpage

\section{Introduction} \label{introduction}
Multi-criteria decision-making methods assist decision-makers in ranking a set of alternatives based on evaluations across multiple criteria. The TODIM method (an acronym in Portuguese for interactive and multiple criteria decision-making) \citep{gl1991}, which generates a linear order of alternatives with ties, ranks alternatives based on the dominance degree of one alternative over another, which is quantified using a dominance function. In the TODIM method, the weights of the criteria are required. The entropy weight method, which derives objective weights from the information entropy of evaluations \citep{shannon1948}, is commonly used in conjunction with the TODIM method \citep{hzl2024}. This paper focuses on the entropy-weight TODIM method. Given that TODIM accounts for psychological biases of decision-makers, such as loss aversion and risk aversion, it has been widely applied in various fields, including waste-to-energy decision-making \citep{mo2022, wu2023, nar2022}, renewable energy investment decision-making \citep{wu2019, has2021}, renewable energy site selection \citep{zhao2024, wang2021}, renewable energy risk assessment \citep{he2022, yin2022, liu2024}, and performance evaluation of renewable energy \citep{zhang2022}.

It has been observed that multi-criteria decision-making methods are susceptible to rank reversal, a phenomenon where the order of two alternatives is altered after the addition or removal of another alternative. This phenomenon was first identified in the analytic hierarchy process \citep{bg1983} and later examined in the TODIM method \citep{jlb2024}. Research on rank reversal has primarily focused on single decision-making methods, such as the preference ranking organization method of enrichment evaluations \citep{ll2024}, the elimination and choice expressing reality \citep{lm2021}, and the technique for order preference by similarity to ideal solution \citep{yzz2022}. \textbf{As far as we know, few studies have addressed how to reduce the rank reversal probability of hybrid methods, such as the entropy-weight TODIM method}.

To address the aforementioned issue, this paper proposes an expanded evaluation matrix to reduce the probability of rank reversal in the entropy-weight TODIM method. To achieve this, virtual alternatives are introduced to construct the expanded evaluation matrix. Subsequently, simulation experiments are conducted to evaluate the performance of the expanded evaluation matrix in mitigating rank reversal.

The contributions and innovations of this paper are summarized as follows.

(1) The expanded evaluation matrix is introduced for the entropy-weight TODIM method to reduce the rank reversal probability. First, equally spaced virtual evaluations are defined, from which virtual alternatives are generated. Next, the expanded evaluation matrix is constructed using these virtual alternatives. When applying the entropy-weight TODIM method, replacing the original evaluation matrix with the expanded evaluation matrix helps reduce the probability of rank reversal. \textbf{To the best of our knowledge, the proposed expanded evaluation matrix is the first method to reduce the probability of rank reversal for hybrid multi-criteria decision-making methods}.

(2) Simulations are conducted to demonstrate the effectiveness of the expanded evaluation matrix in mitigating rank reversal. This paper examines four simulation factors: the number of alternatives, the number of criteria, the level parameter, and the number of added alternatives. Based on the simulation results, several findings and recommendations are provided for the application of the expanded evaluation matrix.

(3) A case study on energy storage technology selection is presented to explore the potential of the proposed expanded evaluation matrix in real-world applications. Sensitivity and comparative analyses are conducted to demonstrate the reliability of the proposed method.

The structure of this paper is organized as follows: 
Section~\ref{section: preliminaries} reviews the research on rank reversal of the TODIM method. Section~\ref{section: expanded evaluation matrix} introduces the extended evaluation matrix, followed by an example of using the extended evaluation matrix. Section~\ref{section: performance} establishes the simulation experiments. Section~\ref{section: case} presents the case study on the selection of energy storage technologies. Section~\ref{section: discussion} provides some suggestions, and the paper closes with conclusions in Section~\ref{section: conslusion}.

\section{Literature review} \label{section: preliminaries}
By using the prospect theory \citep{kt1979,tk1992} to deal with risk aversion and the loss aversion of decision-makers, \citet{gl1991} proposed a multi-criteria decision-making method named TODIM. To apply the TODIM method, the evaluations under each criterion are standardized through standardization techniques. Then, based on the standardized evaluations, a preference function is applied to calculate the dominance degree of each pair of alternatives on every criterion. After that, the total dominance degree of every pair of alternatives is calculated by taking a weighted sum of the dominance degrees on all criteria. Based on the total dominance degrees, the scores of alternatives are calculated, which are used to rank alternatives in descending order. The procedure of the TODIM method is given in \nameref{Appendix}. Since the TODIM method was introduced, it has been widely applied in many fields like scheduling \citep{hdl2024}, emergency management \citep{ctz2024,ltz2024}, energy management \citep{lml2024,zlz2024}, and global climate management \citep{zwc2024}.

To perform the weighted sum in the TODIM method, the weights of criteria are required. In addition to allowing decision-makers to directly give weights, the entropy weight method is another way to generate weights \citep{hzl2024}. The basis of entropy weight is information entropy \citep{shannon1948}. For a criterion, the greater the entropy of the evaluations under the criterion, the more information the evaluations can provide, and the more important the criterion is. The procedure of the entropy weight method is summarized in \nameref{Appendix}. 

As research on the entropy-weight TODIM method continues to advance, \citet{jlb2024} found that when a new alternative is considered when using the entropy-weight TODIM method, the orders of some pairs of original alternatives might be changed, a phenomenon known as rank reversal. Referring to \citet{jlb2024}, the definition of rank reversal is given as follows.
\begin{definition} \label{definition: rank reversal}
\em
    \textbf{\em (Rank reversal).} When applying the entropy-weight TODIM method to a set of alternatives $A=\{a_1,a_2,\ldots,a_n\}$, the strict preference relation is the antisymmetric relation $\textbf{P}_A = \{(a_{i_1}, a_{i_2}) \in A^2 \mid a_{i_1} \mbox{ is strictly preferred to } a_{i_2}\}$; the indifference relation is the symmetric relation $\textbf{I}_A = \{(a_{i_1}, a_{i_2}) \in A^2 \mid a_{i_1} \mbox{ and } a_{i_2} \mbox{ are indifferent}\}$. Let $\hat A$ be the set of alternatives after adding an alternative to $A$. Rank reversal occurs if there exists a pair of alternatives $a_{i_1}$ and $a_{i_2}$ in $A \cap {\hat A}$ satisfying one of the following conditions:

    1) $a_{i_1}, a_{i_2} \in \textbf{P}_A$ while $a_{i_2}, a_{i_1} \in \textbf{P}_{\hat A}$;
    
    2) $a_{i_1}, a_{i_2} \in \textbf{P}_A$ while $a_{i_1}, a_{i_2} \in \textbf{I}_{\hat A}$;
    
    3) $a_{i_1}, a_{i_2} \in \textbf{I}_A$ while $a_{i_2}, a_{i_1} \in \textbf{P}_{\hat A}$;
    
    4) $a_{i_1},a_{i_2} \in \textbf{I}_A$ while $a_{i_1},a_{i_2} \in \textbf{P}_{\hat A}$.
\end{definition}

In the field of multi-criteria decision-making, the study on rank reversal mainly focuses on single decision-making methods like the preference ranking organization method of enrichment evaluations \citep{ll2024}, the elimination and choice expressing reality \citep{lm2021}, and the technique for order preference by similarity to ideal solution \citep{yzz2022}. For the hybrid methods like the entropy-weight TODIM method, \citet{jlb2024} performed a simulation study to study the rank reversal probability of the entropy-weight TODIM method, and investigated the factors affecting the rank reversal probability. However, \citet{jlb2024} did not propose methods to reduce the rank reversal probability of the entropy-weight TODIM method. \textbf{As far as we know, how to reduce the rank reversal probability of the hybrid methods like the entropy-weight TODIM method is still a problem to be addressed.}

\section{Expanded evaluation matrix} \label{section: expanded evaluation matrix}
This section proposes the expanded evaluation matrix to reduce the rank reversal probability of the entropy-weight TODIM method. In Section~\ref{section: decision problem}, the decision problem is first described. The method to construct the expanded evaluation matrix is introduced in Section~\ref{section: structure of the expanded matrix}, followed by an example to demonstrate the effectiveness of the proposed method.

\subsection{Description of the decision problem} \label{section: decision problem}
This paper explores the problem that a decision-maker needs to choose the optimal alternative from $n$ alternatives denoted by $A=\{a_i \mid i=1,2,\ldots,n\}$. The known information is the evaluations of alternatives on $m$ criteria represented by $C=\{c_j \mid j=1,2,\ldots,m\}$. In this paper, $x^j_i \in [0,1]$ represents the evaluation of alternative $a_i$ on criterion $c_j$. The original evaluation matrix is written as follows.
\begin{equation*}
    \textit{X}=
    \begin{matrix}
     & \begin{matrix} c_1 & c_2 & \ldots & c_m \end{matrix} \\
    \begin{matrix} a_1 \\ a_2 \\ \vdots \\ a_n \end{matrix}
    & \begin{bmatrix}
        x^1_1          & x^2_1          & \ldots & x^2_1 \\
        x^1_2          & x^2_2          & \ldots & x^2_2 \\
        \vdots         & \vdots         & \ddots & \vdots \\
        x^1_n          & x^2_n          & \ldots & x^2_n \\
    \end{bmatrix}
\end{matrix}
\end{equation*}

In multi-criteria decision-making, if the evaluations on a criterion are all the same, no useful information is provided on this criterion. To avoid such meaningless situations, this paper requires the evaluations on each criterion to have differences, which means that the maximum evaluation is different from the minimum evaluation on each criterion.

\subsection{Construction of the expanded evaluation matrix} \label{section: structure of the expanded matrix}
Let $x^j_{max}$ (resp.\ $x^j_{min}$) be the maximum (resp.\ minimum) evaluation on criterion $c_j$. For criterion~$c_j$, it is easy to find $v$ equally spaced virtual evaluations, denoted by $x^j_{d,v}$ ($d=1,2,\ldots,v$), between $x^j_{max}$ and $x^j_{min}$. The formula to calculate the equally spaced virtual evaluations is given as
\begin{equation} \label{eq: virtual evaluations}
    x^j_{d,v} = x^j_{min} + \frac{d-1}{v-1} \times (x^j_{max} - x^j_{min})
\end{equation}
In Eq~(\ref{eq: virtual evaluations}), $d$ is the position parameter. The equally spaced virtual evaluation $x^j_{d,v}$ rises when the position parameter $d$ increases. 

For the same position parameter $d$, there are $m$ equally spaced virtual evaluations $x^j_{d,v}$ ($j=1,2,\ldots,m$) corresponding to $m$ criteria, which can be seen as $m$ evaluations of a virtual alternative $\hat a_{d,v}$. Given that $d=1,2,\ldots,v$, there are $v$ virtual alternatives, which are put into a virtual alternative set $\hat A_v = \{\hat a_{d,v} \mid d=1,2,\ldots,v \}$. In this paper, $v$ is the level parameter of the virtual alternative set $\hat A_v$. By putting the virtual alternative set $\hat A_v$ to the original alternative set $A$, the expanded alternative set is set up, denoted by $\bar A_v = A \cup \hat A_v$. The corresponding expanded evaluation matrix of $\bar A_v$ is shown as follows
\begin{equation*}
    \bar{\textit{X}}_v=
    \begin{matrix}
     & \begin{matrix} c_1 & ~\,c_2 & \ldots & ~~c_m \end{matrix} \\
    \begin{matrix} a_1 \\ a_2 \\ \vdots \\ a_n \\ \hat a_{1,v} \\ \hat a_{2,v}\\ \vdots \\ \hat a_{v,v} \end{matrix}
    & \begin{bmatrix}
        x^1_1          & x^2_1          & \ldots & x^2_1 \\
        x^1_2          & x^2_2          & \ldots & x^2_2 \\
        \vdots         & \vdots         & \ddots & \vdots \\
        x^1_n          & x^2_n          & \ldots & x^2_n \\
        x^1_{1,v}      & x^2_{1,v}      & \ldots & x^m_{1,v} \\
        x^1_{2,v}      & x^2_{2,v}      & \ldots & x^m_{2,v} \\
        \vdots         & \vdots         & \ddots & \vdots \\
        x^1_{v,v}      & x^2_{v,v}      & \ldots & x^m_{v,v} \\
    \end{bmatrix}
\end{matrix}
\end{equation*}

\begin{example} \label{example: expanded evaluation matrix}
\em
    \textit{\textbf{(Example of a $2$-level expanded evaluation matrix).}} An evaluation matrix for two alternatives and three criteria is shown as
    \begin{equation*}
    \textit{X}=~~~
    \begin{matrix}
     & \begin{matrix} c_1 & ~c_2 & ~c_3 \end{matrix} \\
    \begin{matrix} a_1 \\ a_2\\ a_3 \end{matrix}
    & \begin{bmatrix}
        0.1147   & 0.9913   & 0.4275 \\
        0.7880   & 0.4822   & 0.2693 \\
        0.3189   & 0.3297   & 0.2803
    \end{bmatrix}
    \end{matrix}\,.
    \end{equation*}

    Referring to Eq~(\ref{eq: virtual evaluations}), for criterion $c_1$, the equally spaced virtual evaluations are calculated as $0.1147+0 \div 1 \times (0.7880-0.1147) = 0.1147$ and $0.1147+1 \div 1 \times (0.7880-0.1147) = 0.7880$. For criterion $c_2$ (resp.\ $c_3$), the equally spaced virtual evaluations are $0.3297$ and $0.9913$ (resp.\ $0.2693$ and $0.4275$). Based on these equally spaced virtual evaluations, the expanded evaluation matrix is constructed as follows.
    \begin{equation*}
        \bar{\textit{X}}_2=~~~
        \begin{matrix}
         & \begin{matrix} c_1 & ~~~c_2~~ & ~~c_3~ \end{matrix} \\
        \begin{matrix} a_1 \\ a_2 \\a_3 \\ \hat a_{1,2} \\ \hat a_{2,2} \end{matrix}
        & \begin{bmatrix}
            0.1147   & 0.9913   & 0.4275 \\
            0.7880   & 0.4822   & 0.2693 \\
            0.3189   & 0.3297   & 0.2803 \\
            0.7880   & 0.9913   & 0.4275 \\
            0.1147   & 0.3297   & 0.2693 \\
        \end{bmatrix}
        \end{matrix}
    \end{equation*}

    For the original evaluation matrix $\textit{X}$, the ranking generated through the entropy-weight TODIM method is $a_1 > a_2 > a_3$, where $a_1 > a_2$ means alternative $a_1$ dominates alternative $a_2$. When another alternative 
    $a_4$ with evaluations $0.8043$, $0.091$, and $0.9674$ is added, the ranking is changed to $a_2 > a_1 > a_3$. Given that the order between $a_1$ and $a_2$ is changed from $a_1 > a_2$ to $a_2 > a_1$, according to Definition~\ref{definition: rank reversal}, the rank reversal phenomena occurs.

    For the expanded evaluation matrix $\bar{\textit{X}}_2$, the ranking is $\hat a_{1,2} > a_1 > a_2 > a_3 > \hat a_{2,2}$. After $a_4$ is added, the ranking is changed to $\hat a_{1,2} > a_4 > a_1 > a_2 > a_3 > \hat a_{2,2}$. Given that the order between $a_1$ and $a_2$ is constant before and after adding $a_4$, according to Definition~\ref{definition: rank reversal}, the rank reversal phenomena does not occur.
\end{example}

Two features are found in Example~\ref{example: expanded evaluation matrix}: (a) adding two virtual alternatives $\hat a_{1,2}$ and $\hat a_{2,2}$ does not change the ranking of $a_1$, $a_2$, and $a_3$; (b) when the expanded evaluation matrix is applied, the order of $a_1$ and $a_2$ is not changed after $a_4$ is added. From these two features, it is easy to put forward the following two hypotheses:
\begin{hypothesis} \label{hyp: consistency}
    \textbf{\textit{(Consistency hypothesis).}} Using the expanded evaluation matrix does not change the ranking of the original alternatives. 
\end{hypothesis}
\begin{hypothesis} \label{hyp: anti-reversal}
    \textbf{\textit{(Anti-reversal hypothesis).}} Rank reversal caused by adding an alternative is avoided when the expanded evaluation matrix is applied. 
\end{hypothesis}

\section{Performance of the expanded evaluation matrix} \label{section: performance}
In Section~\ref{section: structure of the expanded matrix}, two hypotheses were refined from Example~\ref{example: expanded evaluation matrix}. The question is whether these two hypotheses hold in any decision-making situation. To answer this question, this section set up simulations to perform further analyses.
\subsection{Design of simulation experiments} \label{section: simulation design}
\subsubsection{Simulation parameters}
This study examines four key simulation parameters: the number of alternatives, the number of criteria, the level parameter in the expanded evaluation matrix, and the number of added alternatives. Specifically, the number of alternatives and the number of criteria are set to $3, 6, 9, 11, 15$; the level parameter is selected from $2, 4, 6, 8, 10, 15, 20$; and the number of added alternatives is set to $1, 3, 5, 7, 9$. The table below presents a summary of the simulation parameter settings.
\begin{table}[htbp]
\caption{Setting of simulation parameters}
\label{table: parameters}
\centering
\small
\begin{tabular}{llll}
\hline
Number of Alternatives & Number of Criteria & Level parameter & Number of Added Alternatives \\ \hline
3, 6, 9, 11, 15 & 3, 6, 9, 11, 15 & 2, 4, 6, 8, 10, 15, 20 & 1, 3, 5, 7, 9 \\ \hline
\end{tabular}
\end{table}

\subsubsection{Simulation steps}
\textbf{Step 1. Create original evaluation matrices.} Given that there are $5$ values for the number of alternatives and $5$ values for the number of criteria, $5 \times 5 = 25$ simulations are set up. For each simulation, $500$ replicates are performed. One evaluation matrix is generated in every replicate, resulting in $500$ original evaluation matrices. For all simulations, $25 \times 500 = 12500$ original evaluation matrices are generated, which are put into the original matrix set $\textbf{I}_o$. In this paper, evaluation values are sampled from a uniform distribution over $[0,1]$, and the evaluations for different criteria are independently generated to ensure no interaction between criteria.

\textbf{Step 2. Compute original rankings.} For each original evaluation matrix, an original ranking is generated through the entropy-weight TODIM method. For all simulations, $25 \times 500 = 12500$ original rankings are obtained, which are put into the original ranking set $\textbf{T}_o$.

\textbf{Step 3. Obtain expanded evaluation matrices.} For every original evaluation matrix, using Eq.~(\ref{eq: virtual evaluations}) to generate the expanded evaluation matrices. Given that there are $7$ values for the level parameter, each original evaluation matrix corresponds to $7$ expanded evaluation matrices. For all simulations, $25 \times 500 \times 7 = 87500$ expanded evaluation matrices are created, which are put into an expanded evaluation matrix set $\textbf{I}_a$.

\textbf{Step 4. Compute expanded rankings.} Employing the entropy-weight TODIM method, an expanded ranking is generated from an expanded evaluation matrix. Given that there are $87500$ expanded evaluation matrices, $25 \times 500 \times 7 = 87500$ expanded rankings are obtained, which are put into an expanded ranking set $\textbf{T}_a$.

\textbf{Step 5. Add random alternatives.} For each original evaluation matrix, added alternatives are randomly generated based on a uniform distribution over $[0,1]$. Given that the number of added alternatives can be $1$, $3$, $5$, $7$, and $9$, one original evaluation matrix corresponds to $5$ after-addition evaluation matrices. For all simulations, there are $25 \times 500 \times 5 = 62500$ after-addition evaluation matrices, which are put into an after-addition matrix set $\textbf{I}^{add}_o$. Similarly, for each expanded evaluation matrix, there are also $5$ corresponding after-addition expanded evaluation matrices. For all simulations, $25 \times 500 \times 5 \times 7 = 437500$ after-addition expanded evaluation matrices are obtained, which are put into an after-addition expanded evaluation matrix set $\textbf{I}^{add}_a$.

\textbf{Step 6. Compute after-addition rankings.} For each after-addition evaluation matrix (resp.\ each after-addition expanded evaluation matrix), an after-addition ranking (resp.\ after-addition expanded ranking) can be generated through the entropy-weight TODIM method. For all simulations, there are $62500$ after-addition rankings (resp.\ $437500$ after-addition expanded rankings), which are put into an after-addition ranking set $\textbf{T}^{add}_o$ (resp.\ an after-addition expanded ranking set $\textbf{T}^{add}_a$).

\textbf{Step 7: Calculate consistency probability} By comparing $\textbf{T}_o$ with $\textbf{T}_a$, the replicates, where the original ranking in $\textbf{T}_o$ does not match the expanded rankings in $\textbf{T}_a$, can be identified. Let $V_c$ be the number of those inconsistent replicates. The consistency probability is determined by:
\begin{equation} \label{eq: consistency probability}
    P_c = \frac{1-V_c}{500}\,.
\end{equation}

For each inconsistent replicate, the pairs of reversed alternatives can be identified through Definition~\ref{definition: rank reversal}. The total number of pairs of reversed alternatives can be calculated for all inconsistent replicates, denoted by $\varphi_c$. The mean of the number of pairs of reversed alternatives $\bar \varphi_c$ (pair-number mean for short) can be calculated by:
\begin{equation} \label{eq: inconsistency degree}
    \bar \varphi_c = \frac{\varphi_c}{500}\,.
\end{equation}

\textbf{Step 8: Calculate rank reversal probability.} By comparing $\textbf{T}_o$ with $\textbf{T}^{add}_o$ (or comparing $\textbf{T}^{add}_o$ with $\textbf{T}^{add}_a$), the replicates in which rank reversal occurs can be identified. Let $V_r$ be the number of those replicates. The rank reversal probability can be computed by:
\begin{equation} \label{eq: rank reversal probability}
    P_r = \frac{V_r}{500}\,.
\end{equation}

\subsection{Experimental results} \label{section: results}
\subsubsection{Disscusion on the consistency probability} \label{section: discuss on consistency}
In this section, the impact of the expanded evaluation matrix on the original ranking is investigated by analyzing the consistency probability.

When studying the consistency probability, the simulation experiment includes $5$ values for the number of alternatives, $5$ values for the number of criteria, and $7$ values for the level parameter. It is easy to find that there are $5 \times 5 \times 7 = 175$ combinations, corresponding to $175$ consistency probabilities and $175$ pair-number means. The average (resp.\ standard variance) of the $175$ inconsistency probabilities is $0.469$. The average (resp.\ standard variance) of the $175$ pair-number means is $0.81$ (resp.\ $0.907$).

It is found that the consistency probabilities are not high, which means using the expanded evaluation matrix might cause changes in the original rankings. Therefore, Hypothesis~\ref{hyp: consistency} is invalid. Although the consistency probabilities are not high, the pair-number mean is low, which means using the expanded evaluation matrix influences a very small number of pairs of alternatives.

\subsubsection{Discussion on the rank reversal probability} \label{section: discuss on reversal probability}
This section investigates how effective the expanded evaluation matrix is in reducing the rank reversal probability of the entropy-weight TODIM method.

When studying the rank reversal probability, the simulation experiment includes $5$ values for the number of alternatives, $5$ values for the number of criteria, $7$ values for the level parameter, and $5$ values for the number of added alternatives, so there are $5 \times 5 \times 7 \times 5= 875$ combinations. For each combination, we can obtain one rank reversal probability with the original evaluation matrix and one rank reversal probability with the expanded evaluation matrix. The former probability minus the latter probability is a probability gap. The results show that all $875$ probability gaps are positive, so Hypothesis~\ref{hyp: anti-reversal} is valid. The average (resp.\ standard variance) of the $875$ probability gaps is $0.1032$ (resp.\ $0.0604$). For different numbers of added alternatives, the average and standard variance of the probability gaps are given in Table~\ref{tab: rank reversal probability}.
\begin{table}[htbp]
\centering
\caption{Summary for different numbers of added alternatives.}
\label{tab: rank reversal probability}
\begin{tabular}{llllll}
\hline
Number of added alternatives & 1 & 3 & 5 & 7 & 9 \\
Average of probability gaps & 0.1159 & 0.118 & 0.1027 & 0.0934 & 0.086 \\
Standard variance of probability gaps & 0.0679 & 0.0658 & 0.0571 & 0.054 & 0.0492 \\
\hline
\end{tabular}
\end{table}

Based on the above analyses, it can be found that the rank reversal probability of the entropy-weight TODIM method is reduced by using the expanded evaluation matrix. One thing worth noting is that even if more than one alternative is added, the expanded evaluation matrix still works (referring to the last column of Table~\ref{tab: rank reversal probability}). Given that the average rank reversal probability of the entropy-weight TODIM method is $0.5311$, the reduction of $0.1032$ ($19.43\%$) is remarkable.

Subsequently, three kinds of correlation tests ({\em i.e.}, Pearson correlation test, Spearman correlation test, and the Kendall correlation test) are performed to explore the relationship between the probability gap and three parameters: the number of alternatives, the level parameter, and the number of added alternatives. The results are shown in Table~\ref{tab: correlation}.
\begin{table}[htbp]
\centering
\caption{Correlation coefficients of the correlation tests.}
\label{tab: correlation}
\begin{tabular}{llll}
\hline
		           & Number of alternatives & Level parameter & Number of added alternatives \\ \hline
		Pearson  & 0.26$^{***}$           & 0.7207$^{***}$  & -0.1978$^{***}$              \\
		Spearman & 0.26$^{***}$           & 0.7601$^{***}$  & -0.1742$^{***}$              \\
		Kendall  & 0.1979$^{***}$         & 0.6044$^{***}$  & -0.1285$^{***}$             \\ \hline
              \multicolumn{4}{l}{Note. ``\,$^{***}$\," denotes the $1\%$ significance level.} \\
	\end{tabular}
\end{table}

From Table~\ref{tab: correlation}, it is found that there are positive correlations between the probability gap and the number of alternatives (resp.\ level parameter), which means when the number of alternatives (resp.\ the number of virtual alternatives) is large, the effect of using the expanded evaluation matrix to reduce the rank reversal probability is obvious. In the last column of Table~\ref{tab: correlation}, the correlation coefficients are all negative, which means there are negative correlations between the probability gap and the number of added alternatives. Hence, when the number of added alternatives is large, the effect of using the expanded evaluation matrix to reduce the rank reversal probability is limited.
\section{Case study} \label{section: case}
This section presents a case study on selecting energy storage technology to showcase the potential of the expanded evaluation matrix in addressing real-world issues. Sensitivity analyses and comparative analyses are performed to illustrate the reliability of the proposed method.

\subsection{Background} \label{section: case background}
Renewable energy sources, such as solar and wind power, are inherently intermittent and subject to fluctuations. Energy storage technologies are employed to store surplus electricity and release it during periods of peak demand, thereby mitigating the temporal mismatch between energy production and consumption and reducing associated energy losses. Due to these advantages, energy storage technologies have been extensively applied in peak shaving, frequency regulation, grid backup power, and the integration of renewable energy \citep{ybm2022,whb2024}.

Given the complexity and variety of application scenarios for energy storage technology, selecting the most suitable energy technology is crucial for each specific context. Considering the multiple criteria involved in selecting energy storage technology, this issue can be framed as a multiple-criteria decision-making problem. Several methods, including the analytical hierarchy process \citep{lkf2021, aml2023}, the technique for order of preference by similarity to ideal solution \citep{tdi2021, ljd2022}, the Dempster–Shafer evidence theory \citep{hwh2020}, and the preference ranking organization method of enrichment evaluations \citep{wzg2021}, have been applied to address this challenge.

Due to the extended construction periods of energy storage projects, new technologies may emerge during construction. For instance, in 2020, a multiple-criteria decision-making method was used to select the most suitable technology from five energy storage alternatives. Subsequently, construction of the project commenced. During the project's construction, another technology emerges. If the new technology is also considered in 2020, would the initially selected technology still be the most suitable? The rank reversal phenomenon inherent in multiple-criteria decision-making methods means that the addition of alternatives can result in a shift in the optimal alternative (as illustrated in Example~\ref{example: expanded evaluation matrix}). To mitigate the impact of adding alternatives on decision outcomes, the probability of rank reversal must be minimized across various fields \citep{sen2021, cck2014}. However, in the context of energy storage technology selection, few studies address the reduction of rank reversal probability when applying multiple-criteria decision-making methods. To bridge this research gap, the entropy-weight TODIM method, coupled with the expanded evaluation matrix, is applied to the energy storage technology selection problem in this section.

\subsection{Data collection and problem solving} \label{section: case solving}
This paper considers ten energy storage technologies: pumped storage ($a_1$), compressed air energy storage ($a_2$), flywheel energy storage ($a_3$), sodium-sulfur battery ($a_4$), vanadium flow battery ($a_5$), colloid battery ($a_6$), lead-carbon battery ($a_7$), lithium iron phosphate battery ($a_8$), superconducting energy storage ($a_{9}$), and supercapacitor ($a_{10}$).

These technologies are assessed based on fifteen criteria: power level ($c_1$), response speed ($c_2$), continuous discharge time ($c_3$), discharge depth ($c_4$), cycle number ($c_5$), energy conversion efficiency ($c_6$), self-discharge rate ($c_7$), volume power density ($c_8$), volume energy density ($c_9$), environmental impact ($c_{10}$), capacity unit price ($c_{11}$), power unit price ($c_{12}$), operation and maintenance cost ($c_{13}$), security ($c_{14}$), and technology maturity ($c_{15}$).

The normalized evaluations of these ten technologies on the fifteen criteria, sourced from \citet{hwh2020}, are presented in Table~\ref{table: evaluation}.
\begin{table}[htbp]
	\centering
        \small
        \caption{The normalized evaluations of eleven types of energy storage technologies.}
	\label{table: evaluation}
	\begin{tabular}{llllllllllllllll}
		\hline
		 &     $c_1$&     $c_2$&          $c_3$&          $c_4$&          $c_5$&          $c_6$&          $c_7$&          $c_8$&          $c_9$&          $c_{10}$&          $c_{11}$&          $c_{12}$&     $c_{13}$&          $c_{14}$&$c_{15}$     \\ \hline
		 $a_1$& 1  & 1  & 1  & 0  & 0.81  & 0.46  & 1  & 0.28  & 1  & 1  & 0  & 0  & 1  & 1  & 1  \\
		 $a_2$& 1  & 1  & 1  & 0  & 0.19  & 0  & 1  & 0.78  & 1  & 1  & 0  & 0.01  & 0  & 0.67  & 0.5  \\
		 $a_3$& 0.5  & 1  & 0  & 0  & 0.81  & 0.73  & 0  & 0.56  & 0.74  & 1  & 0.28  & 0.16  & 1  & 0.33  & 0.5  \\
		 $a_4$& 0.5  & 0.5  & 0.67  & 1  & 0  & 0.66  & 1  & 0.89  & 0.99  & 1  & 0.01  & 0.75  & 1  & 0  & 1  \\
		 $a_5$& 0.5  & 0.5  & 1  & 1  & 0.22  & 0.05  & 0.67  & 0  & 0.99  & 1  & 0  & 0.06  & 0.5  & 0.67  & 0.5  \\
		 $a_6$& 0.5  & 0.5  & 1  & 0.33  & 0.03  & 0.39  & 1  & 1  & 1  & 1  & 0.02  & 0.25  & 0.5  & 0.33  & 1  \\
		 $a_7$& 0.5  & 0.5  & 0.67  & 0.33  & 0.01  & 0.8  & 1  & 1  & 1  & 1  & 0.02  & 0.13  & 0.5  & 0.67  & 0  \\
		 $a_8$& 0.5  & 0.5  & 0.67  & 0.67  & 0.05  & 0.93  & 0.67  & 1  & 1  & 1  & 0.33  & 1  & 0.5  & 0.33  & 1  \\
		 $a_{9}$& 0.5  & 0  & 0  & 0  & 0.19  & 0.53  & 0.67  & 0.42  & 0  & 1  & 0.15  & 0  & 0  & 0.33  & 0.5  \\
		 $a_{10}$& 0  & 0  & 0.33  & 0  & 1  & 0.80  & 0.33  & 1  & 0.97  & 0  & 0.32  & 0  & 0.5  & 0.33  & 0.5  \\ \hline
	\end{tabular}
\end{table}

The entropy-weight TODIM method, along with the expanded evaluation matrix, is employed in this analysis. The level parameter of the expanded evaluation matrix is set to $8$. As outlined in Eq.~(\ref{eq: virtual evaluations}), eight virtual alternatives are generated. Utilizing the entropy-weight TODIM method, the ranking of the ten energy storage technologies is determined as follows: $R_1 = a_8>a_1>a_4>a_6>a_7>a_2>a_3>a_5>a_{10}>a_9$. Based on this ranking, the lithium iron phosphate battery $a_8$ is identified as the optimal choice.

Suppose that during the construction of the lithium iron phosphate battery project, two new technologies are introduced, with the evaluations given by the following sets: $\{0.02$, $0.81$, $0.66$, $0.54$, $0.81$, $0.36$, $0.85$, $0.4$, $0.58$, $0.83$, $0.04$, $0.37$, $0.76$, $0.96$, $0.23\}$ and $\{0.63$, $0.81$, $0.9$, $0.26$, $0.53$, $0.1$, $0.96$, $0.61$, $0.82$, $0.5$, $0.96$, $0.19$, $0.69$, $0.9$, $0.09\}$. If these new technologies are included in the analysis, the after-addition ranking of the original ten technologies remains $R_2 = a_8>a_1>a_4>a_6>a_7>a_2>a_3>a_5>a_{10}>a_9$, as determined by the entropy-weight TODIM method with the expanded evaluation matrix. It can be found that $R_1$ and $R_2$ are consistent and $a_8$ ranks first in both rankings. Therefore, it is reliable to select the lithium iron phosphate battery as the optimal choice.

\subsection{Sensitivity analysis} \label{section: case sensitive}
In the previous section, the level parameter was set to $8$. This section investigates the impact of varying the level parameter on the ranking.

First, the level parameter is set to $6$, $7$, $8$, $9$, and $10$. The entropy-weight TODIM method, applied with the expanded evaluation matrix, is then used to obtain the rankings corresponding to each level parameter. The rankings of the original ten technologies, both before and after the inclusion of the two new technologies, are presented in Table~\ref{tab: sensitive}. 
\begin{table}[htbp]
	\centering
        \footnotesize
        \caption{Results of the sensitivity analyses.}
	\label{tab: sensitive}
	\begin{tabular}{lll}
		\hline
	LP & Ranking before adding two new technologies & Ranking after adding two new technologies \\ \hline
        6 & $a_8>a_1>a_4>a_6>a_7>a_2>\mathbf{a_3>a_5}>a_{10}>a_9$ & $a_8>a_1>a_4>a_6>a_7>a_2>\mathbf{a_5>a_3}>a_{10}>a_9$ \\
        7 & $a_8>a_1>a_4>a_6>a_7>a_2>\mathbf{a_3>a_5}>a_{10}>a_9$ & $a_8>a_1>a_4>a_6>a_7>a_2>\mathbf{a_5>a_3}>a_{10}>a_9$ \\
        8 & $a_8>a_1>a_4>a_6>a_7>a_2>a_3>a_5>a_{10}>a_9$ & $a_8>a_1>a_4>a_6>a_7>a_2>a_3>a_5>a_{10}>a_9$ \\
        9 & $a_8>a_1>a_4>a_6>a_7>a_2>a_3>a_5>a_{10}>a_9$ & $a_8>a_1>a_4>a_6>a_7>a_2>a_3>a_5>a_{10}>a_9$ \\
        10 & $a_8>a_1>a_4>a_6>a_7>\mathbf{a_3>a_2}>a_5>a_{10}>a_9$ & $a_8>a_1>a_4>a_6>a_7>\mathbf{a_2>a_3}>a_5>a_{10}>a_9$ \\ \hline
        \multicolumn{3}{l}{Note. ``LP" represents the level parameter.} \\
	\end{tabular}
\end{table}

The differences between the rankings before and after the addition of the new technologies are highlighted in Table~\ref{tab: sensitive}. It can be found that when the level parameter is $6$ and $7$ (resp.\ $10$), the order of $a_3$ and $a_5$ (resp.\ the order of $a_2$ and $a_3$) is reversed after the inclusion of the two new technologies. Notably, $a_8$ consistently ranks as the best alternative, irrespective of the level parameter's value.

\subsection{Comparative analysis} \label{section: case compare}
In this section, we recalculate the ranking of the original ten technologies without using the expanded evaluation matrix. Before the inclusion of the two new technologies, the generated ranking is $R_3 = a_8>a_1>a_4>a_6>a_7>a_2>a_5>a_3>a_{10}>a_9$. A comparison of $R_1$ and $R_3$ reveals that the relative order of $a_2$ and $a_5$ differs between the two rankings; however, $a_8$ consistently ranks as the best alternative in both two rankings. When the two new technologies are included in the analysis, the ranking shifts to $R_4 = a_1>a_8>a_4>a_6>a_7>a_2>a_5>a_3>a_{10}>a_9$. 

It has been observed that, without the expanded evaluation matrix, the entropy-weight TODIM method generates inconsistent rankings for the best alternative before and after the inclusion of the two new technologies. This inconsistency presents a challenge for decision-makers when selecting the optimal alternative. In contrast, the entropy-weight TODIM method, when applied with the expanded evaluation matrix, consistently identifies $a_8$ as the best alternative. This consistency provides decision-makers with accurate and reliable recommendations for making informed decisions.

\section{Discussion} \label{section: discussion}
The simulation experiment in Section~\ref{section: performance} and the case study in Section~\ref{section: case} demonstrate that using the expanded evaluation matrix, rather than the original evaluation matrix, reduces the probability of rank reversal. As outlined in Eq.~(\ref{eq: virtual evaluations}), the method for generating the virtual alternatives in the expanded evaluation matrix is straightforward and easy to implement. Given the simplicity of transitioning from the original evaluation matrix to the expanded evaluation matrix, it is recommended to integrate the expanded evaluation matrix into the standardization process when applying the entropy-weight TODIM method.

To construct the expanded evaluation matrix, the level parameter must be determined in advance. According to the simulation results in Section~\ref{section: discuss on reversal probability}, a larger level parameter corresponds to a smaller rank reversal probability. Therefore, it is recommended that the level parameter should be maximized when employing the expanded evaluation matrix.

\section{Conclusion} \label{section: conslusion}
This paper introduced the expanded evaluation matrix for the entropy-weight TODIM method. Specifically, equally spaced virtual evaluations were defined for each criterion, from which virtual alternatives were generated. By adding these virtual alternatives to the original set of alternatives, the expanded evaluation matrix was established. Simulation experiments revealed that although the rankings of alternatives before and after the inclusion of virtual alternatives may be inconsistent, the number of pairs of reversed alternatives remains small. The simulation results indicated that the expanded evaluation matrix effectively reduced the probability of rank reversal. A case study on selecting energy storage technology demonstrated that the entropy-weight TODIM method, when combined with the expanded evaluation matrix, produced a more reliable ranking of alternatives compared to the entropy-weight TODIM method without the expanded evaluation matrix. Regarding the determination of the level parameter for the expanded evaluation matrix, it is recommended that the level parameter should be as large as possible.

Although this paper introduces the expanded evaluation matrix specifically for the entropy-weight TODIM method, the process of constructing the expanded evaluation matrix is independent of the entropy-weight TODIM method. Therefore, the expanded evaluation matrix might be easily adapted to other hybrid methods. The performance of the expanded evaluation matrix in other hybrid methods can be verified in future research.

\section*{Acknowledgments}
The work was supported by the Project of Sichuan System Science and Enterprise Development Research Center (Xq24B03), and the Energy and Environment Carbon Neutrality Innovation Research Center (YB03202408).

\section*{Disclosure of interest}
We confirm that there are no known conflicts of interest associated with this publication and there has been no significant financial support for this work that could have influenced its outcome.

\section*{Appendix} \label{Appendix}
In a decision problem, a decision maker needs to select an alternative from a set of $n$ alternatives denoted by $A=\{a_i \mid i=1,2,\ldots,n\}$. All $n$ alternatives are evaluated on a set of $m$ criteria denoted by $C=\{c_j \mid i=1,2,\ldots,m\}$. The evaluation of $a_i$ on $c_j$ is denoted by $x^j_i$. Next, the procedure of the entropy-weight TODIM method is recalled step by step.
\begin{procedure}
\textbf{(Procedure of the entropy-weight TODIM method).}

\textbf{Step 1. Normalize the evaluations on each criterion.} The normalized evaluation of evaluation $x^j_i$ is calculated by
\begin{equation*}
    \bar x^j_i = \frac{x^j_i}{\max_{i_1=1,2,\ldots,n}\{x^j_{i_1}\}}\,. 
\end{equation*}

\textbf{Step 2. Calculate the entropy of every criterion.} Let $\hat e_j$ be the information entropy of the evaluations on criterion $c_j$, calculated by
\begin{equation*}
    \hat e_j = \frac{\sum_{i=1,2,\ldots,n}[\bar x^j_i \times ln(\bar x^j_i)]}{ln(n)}
\end{equation*}
where ``$ln$" represnts the natural logarithm.

\textbf{Step 3. Calculate the entropy weights.} The weight of criterion $c_j$ is calculated by
\begin{equation*}
    w_j = \frac{1-\hat e_j}{\sum_{j_1=1,2,\ldots,m}(1 - \hat e_{j_1})}\,.
\end{equation*}

\textbf{Step 4. Calculate the dominance degrees.} The dominance degree of $a_{i_1}$ to $a_{i_2}$ on criterion $c_j$ is calculated by the dominance function
\begin{equation*}
    \hat d^j_{i_1,i_2} = \left\{
        \begin{array}{ll}
            [w_j \times (\bar x^j_{i_1} - \bar x^j_{i_2})]^{0.88} & ,\,if\, \bar x^j_{i_1} - \bar x^j_{i_2} > 0 \\
            0 & ,\,if\, \bar x^j_{i_1} - \bar x^j_{i_2} = 0 \\
            -2.25 \times w_j \times (\bar x^j_{i_2} - \bar x^j_{i_1})^{0.88} & ,\,if\, \bar x^j_{i_1} - \bar x^j_{i_2} < 0 \\
        \end{array}
        \right.
\end{equation*}

\textbf{Step 5. Calculate the total dominance degrees.} The total dominance degree of $a_{i_1}$ to $a_{i_2}$ is calculated by $\bar D_{i_1,i_2} = \sum_{j=1,2,\ldots,m}\hat d^j_{i_1,i_2}$.

\textbf{Step 6. Calculate the score of each alternative.} The score of $a_i$ is calculated by $\hat s_i = \sum_{i_1 = 1,2,\ldots,n} \bar D_{i,i_1}$.

\textbf{Step 7. Alternatives are ranked in the descending order of scores.}
\end{procedure}
\end{document}